\documentclass{amsart}

\usepackage{tikz,amsmath,amstext,amscd,epsfig,euscript, mathrsfs, dsfont,pspicture,multicol,graphpap,graphics,graphicx,times,enumerate,subfig,sidecap,wrapfig,color}
\usepackage{amssymb}

\usepackage[english,francais]{babel}

\newtheorem{theorem}{Theorem}[section]
\newtheorem{lemma}[theorem]{Lemma}
\newtheorem{e-proposition}[theorem]{Proposition}

\newtheorem{e-definition}[theorem]{Definition\rm}

\newcommand{\M}{{\mathcal M}}       %
\newcommand{\R}{{\mathbb R}}       % Field of real numbers
\newcommand{\DD}{{\mathcal D}}

\newcommand{\HH}{{\mathcal H}}

\newcommand{\RR}{{\mathcal R}}

\newcommand{\EE}{{\mathcal E}}

\newcommand{\diam}{\mathop{\rm diam}}

\newcommand{\rf}[1]{{(\ref{#1})}}

\newcommand{\supp}{\operatorname{supp}}

\newcommand{\ve}{{\varepsilon}}

\newcommand{\pom}{{\partial \Omega}}

\newcommand{\LM}{{\mathsf{LM}}}

\def\XXint#1#2#3{{\setbox0=\hbox{$#1{#2#3}{\int}$ }
\vcenter{\hbox{$#2#3$ }}\kern-.58\wd0}}

\begin{document}

\selectlanguage{english}
\title[Rectifiability of harmonic measure]{Harmonic measure is rectifiable if it is absolutely continuous with   respect to the   co-dimension one Hausdorff measure}

\thanks{We acknowledge the support of the following grants: the second author  -- NSF DMS 1361701, the third author -- SEV-2011-0087 and (FP7/2007-2013)/ ERC agreement no. 615112 HAPDEGMT, the fourth author -- Sloan Fellowship, NSF DMS 1344235, NSF DMS 1220089,  NSF  DMR 0212302, the sixth author -- ERC 320501 FP7/2007-2013 (which also funded the first and fifth authors),  2014-SGR-75, MTM2013-44304-P, ITN MAnET (FP7-607647), the last author -- NSF DMS-1265549.
The results of this paper were obtained at the \textit{Institut Henri Poincar\'e}, and at the 2015 ICMAT program \textit{Analysis and Geometry in Metric Spaces}. All authors would like to express their gratitude to these institutions for their support and nice working environments.}

\date{\today}

\selectlanguage{english}
\author[Azzam]{Jonas Azzam}
\address{Jonas Azzam
\\
Departament de Matem\`atiques
\\
Universitat Aut\`onoma de Barcelona
\\
Edifici C Facultat de Ci\`encies
\\
08193 Bellaterra (Barcelona), Catalonia
}
\email{jazzam@mat.uab.cat}

\author[Hofmann]{Steve Hofmann}

\address{Steve Hofmann
\\
Department of Mathematics
\\
University of Missouri
\\
Columbia, MO 65211, USA} \email{hofmanns@missouri.edu}

\author[Martell]{Jos\'e Mar{\'\i}a Martell}

\address{Jos\'e Mar{\'\i}a Martell
\\
Instituto de Ciencias Matem\'aticas CSIC-UAM-UC3M-UCM
\\
Consejo Superior de Investigaciones Cient{\'\i}ficas
\\
C/ Nicol\'as Cabrera, 13-15
\\
E-28049 Madrid, Spain} \email{chema.martell@icmat.es}

\author[Mayboroda]{Svitlana Mayboroda}

\address{Svitlana Mayboroda
\\
Department of Mathematics
\\
University of Minnesota
\\
Minnea\-po\-lis, MN 55455, USA} \email{svitlana@math.umn.edu}

\author[Mourgoglou]{Mihalis Mourgoglou}

\address{Mihalis Mourgoglou
\\
Departament de Matem\`atiques\\
 Universitat Aut\`onoma de Barce\-lo\-na and Centre de Recerca Matem\` atica
\\
Edifici C Facultat de Ci\`encies
\\
08193 Bellaterra (Barcelona), Catalonia
}
\email{mmourgoglou@crm.cat}

\author[Tolsa]{Xavier Tolsa}
\address{Xavier Tolsa
\\
ICREA and Departament de Matem\`atiques
\\
Universitat Aut\`onoma de Bar\-ce\-l\-ona
\\
Edifici C Facultat de Ci\`encies
\\
08193 Bellaterra (Barcelona), Catalonia
}
\email{xtolsa@mat.uab.cat}

\author[Volberg]{Alexander Volberg}
\address{Alexander Volberg
\\
Department of Mathematics
\\
Michigan State University
\\
East Lan\-sing, MI 48824, USA}

\email{volberg@math.msu.edu}

\medskip

\begin{abstract}
\selectlanguage{english}
In the present paper we sketch the proof of the fact that for any open connected set $\Omega\subset\R^{n+1}$, $n\geq 1$, and any $E\subset \partial \Omega$ with $0<\HH^n(E)<\infty$, absolute continuity of the harmonic measure $\omega$ with respect to the Hausdorff measure on $E$ implies that $\omega|_E$ is rectifiable.

%\vskip 0.5\baselineskip

\end{abstract}

\selectlanguage{english}

\maketitle
% main text
\section{Introduction and notation}
\label{intro}

In what follows, $\Omega$ will denote a connected open set.  As usual, using the Perron method and Riesz Representation Theorem, we can define the harmonic measure of $\Omega$ with a pole at $p\in \Omega$ (denote it by $\omega^p$) for any bounded domain and use a limiting procedure to extend this definition to unbounded domains (see, e.g., \cite[Section 3]{HM12}).   

We call a Radon measure $\sigma$ on $\R^{n+1}$  $n$-{\it rectifiable} if its (any) Borel support can be covered by countably many (rotated) graphs of scalar Lipschitz functions on $\R^{n}$ up to zero $\sigma$ measure (in particular any  pure point measure is rectifiable in this definition, but we will apply it to harmonic measures, such measures never have point masses.)

The main result that we announce here is the following
\begin{theorem}
\label{mmain}
Let $\Omega$ be any   open connected set   in $\R^{n+1}$, $n\ge 1$ and let $\omega=\omega^p$, for some $p\in \Omega$, be the harmonic measure on $\Omega$. Let $E$ be a Borel set, and let $\omega|_E\ll\HH^n|_E$.  Then $\omega|_E$ is $n$-rectifiable.
\end{theorem}

A simple corollary is that if $\omega|_E$ and $\HH^n|_E$ are mutually absolutely continuous then $\HH^n|_E$ is $n$-rectifiable. In Besicovitch theory such sets are called $n$-rectifiable.

The structure of harmonic measures was a focus of attention of many mathematicians starting with Carleson's article \cite{C}, which proved that
harmonic measure on any continuum in the plane has dimension  at least $1/2+\epsilon$ for universal positive $\epsilon$. Without an $\epsilon$ the result was known for a long time and was nothing but a re-interpretation of K\"obe distortion theorem. Carleson's result had a very difficult proof, but Makarov found a novel approach that also established that harmonic measure on any continuum on the plane has dimension   exactly  $1$.  Carleson's and Makarov's works generated a fantastic amount of very difficult and revealing results about the structure of harmonic measure.  Topological restrictions (like being a continuum) and dimension played an important part from the start. Let us mention that 1) for any planar domain any harmonic measure has dimension at most $1$
(Jones--Wolff, \cite{JW}); 2) there are no analogs of Makarov's or Jones--Wolff's results in dimension $3$ and higher (Bourgain, \cite{Bo}, Wolff, \cite{W}); 3)
in dimension $2$ (planar domains) if harmonic measure on a continuum dwells on a good set (meaning that its Borel support can be covered by countably many (rotated) Lipschitz graphs) then harmonic measure is absolutely continuous with respect to the {\it surface measure} $\HH^1$ (Bishop--Jones, \cite{BJ}). It is also shown in \cite{BJ} that some topological restrictions are necessary, for there is a 1-rectifiable set in a plane for which 
harmonic measure is positive and Hausdorff measure is zero.

Notice that Theorem \ref{mmain} establishes the converse to Bishop--Jones result: they  claim that in $\R^2$ and with topological restriction (continuum) $1$-rectifiable harmonic measure must be absolutely continuous with respect to  measure $\HH^1$. We claim that without any restrictions whatsoever and in any $\R^{n+1}$ the converse is true: if harmonic measure is absolutely continuous with respect to  $\HH^{n}$, then  it is $n$-rectifiable. Moreover, our result is very local, one can start with harmonic measure restricted to an arbitrary Borel set $E$, and again its absolute continuity with respect to $\HH^{n}|E$ implies that harmonic measure restricted to $E$ is $n$-rectifiable.

In this respect, let us recall that two natural directions were considered in \cite{HM12}, \cite{HMIV}, \cite{HMU}: 1) geometry-to-analysis, when some sort of (say, scalar invariant) rectifiability of the boundary  implies quantitative claims on the density of harmonic measure with respect to $\HH^n$, 2) analysis-to-geometry,  when the absolute continuity of harmonic measure and $\HH^n$ in a scale invariant fashion implies the quantitative rectifiability of the boundary. As mentioned above, in geometry-to-analysis direction  Bishop--Jones in \cite{BJ}  found that  certain additional topological properties  are necessary. But Theorem \ref{mmain} shows that no such obstacles exist in the direction analysis-to-geometry.
 
 This paper is a short announcement of the results in \cite{ahm3tv}, which in turn arises from the union of two separate works, \cite{AMT0} and \cite{HMMTV}. In \cite{AMT0}, a version of Theorem \ref{mmain} was proved under the additional assumption that the boundary of $\Omega$ is {\it porous} in $E$ (a certain topological restriction). In \cite{HMMTV} the porosity assumption was removed. Both   \cite{AMT0} and \cite{HMMTV} exist only as preprints on ArXiv.

 We finish this section with some notation and preliminaries.  Given a signed Radon measure $\nu$ in $\R^{n+1}$ we consider the $n$-dimensional Riesz
transform
$$\RR\nu(x) = \int \frac{x-y}{|x-y|^{n+1}}\,d\nu(y),$$
whenever the integral makes sense. For $\ve>0$, its $\ve$-truncated version is given by 
$$\RR_\ve \nu(x) = \int_{|x-y|>\ve} \frac{x-y}{|x-y|^{n+1}}\,d\nu(y),\,\,\,
R_{*,\delta} \nu(x)= \sup_{\ve>\delta} |\RR_\ve \nu(x)|.
$$

We also consider the maximal operator
$$\M^n_\delta\nu(x) = \sup_{r>\delta}\frac{|\nu|(B(x,r))}{r^n},$$
In the case $\delta=0$ we write $\RR_{*} \nu(x):= \RR_{*,0} \nu(x)$ and $\M^n\nu(x):=\M^n_0\nu(x)$.\\
%If $\mu$ is a fixed positive Radon measure and $f\in L^1_{loc}(\mu)$, the centered maximal Hardy-Littlewood operator applied to $f$ is 
%$$M_\mu f(x) = \sup_{r>0}\frac1{\mu(B(x,r))}\int_{B(x,r)}|f|\,d\mu.$$

Next, let $\HH_\infty^s$ denote the Hausdorff content of order $s$.
We recall a result of Bourgain from \cite{Bo}. 

\begin{lemma}[Bourgain's Lemma]
\label{lembourgain}
There is $\delta_{0}>0$ depending only on $n\geq 1$ so that the following holds for $\delta\in (0,\delta_{0})$. Let $\Omega\subsetneq \R^{n+1}$ be a  bounded domain, $n+1\geq s>n-1$, $\xi \in \partial \Omega$, $r>0$, and $B=B(\xi,r)$. Then
\[ \omega^{x}(B)\gtrsim_{n,s} \frac{\mathcal H_\infty^{s}(\partial\Omega\cap \delta B)}{(\delta r)^{s}}\quad \mbox{  for all }x\in \delta B\cap \Omega .\]
\end{lemma}

  Let us point out that it is not difficult to prove that if Theorem~\ref{mmain} holds on bounded domains, then it remains valid for any open set $\Omega\subset \R^{n+1}$. We shall restrict our discussion here to the case when $\Omega$ is bounded. Also, for brevity, we consider the case $n\geq 2$ only. For details on these and other issues we refer the reader to \cite{ahm3tv}. 

\section{David--Mattila cells, doubling cells}
\label{cells}

What follows is the sketch of the proof of our main result. The details should be found in \cite{ahm3tv}. To start, we  consider the dyadic lattice of cubes
with small boundaries of David-Mattila associated with $\omega^p$. It has been constructed in \cite[Theorem 3.2]{David-Mattila} (with $\omega^p$ replaced by a general Radon measure). 

\begin{lemma}[David, Mattila \cite{David-Mattila}]
\label{lemcubs}
Consider two constants $C_0>1$ and $A_0>5000\,C_0$ and denote $W=\supp\omega^p$. Then there exists a sequence of partitions of $W$ into
Borel subsets $Q$, $Q\in \DD_k$, with the following properties:
%\begin{itemize}
1) For each integer $k\geq0$, $W$ is the disjoint union of the ``cubes'' $Q$, $Q\in\DD_k$, and
if $k<l$, $Q\in\DD_l$, and $R\in\DD_k$, then either $Q\cap R=\varnothing$ or else $Q\subset R$.
 2) For each $k\geq0$ and each cube $Q\in\DD_k$, there is a ball $B(Q)=B(z_Q,r(Q))$ such that
$$z_Q\in W, \quad A_0^{-k}\leq r(Q)\leq C_0\,A_0^{-k},\qquad W\cap B(Q)\subset Q\subset W\cap 28\,B(Q)=W \cap B(z_Q,28r(Q)),$$
and
$\mbox{the balls\, $5B(Q)$, $Q\in\DD_k$, are disjoint.}$
3)   The cubes $Q\in\DD_k$ have ``small boundaries" -- see \cite{David-Mattila} for a precise definition.  
4) Denote by $\DD_k^{db}$ the family of cubes $Q\in\DD_k$ for which
\begin{equation}\label{eqdob22}
\omega^p(100B(Q))\leq C_0\,\omega^p(B(Q)).
\end{equation}
We have that $r(Q)=A_0^{-k}$ when $Q\in\DD_k\setminus \DD_k^{db}$
and 
$\omega^p(100B(Q))\leq C_0^{-l}\,\omega^p(100^{l+1}B(Q))$ for all $l\geq1$ such that $100^l\leq C_0$ and $Q\in\DD_k\setminus \DD_k^{db}$.
\end{lemma}

We use the notation $\DD=\bigcup_{k\geq0}\DD_k$. Observe that the families $\DD_k$ are only defined for $k\geq0$. So the diameters of the cubes from $\DD$ are uniformly
bounded from above.
%For $Q\in\DD$, we set $\DD(Q) =\{P\in\DD:P\subset Q\}$.
  Given $Q\in\DD_k$, we denote $J(Q)=k$.  
We call
$\ell(Q)= 56\,C_0\,A_0^{-k}$ the side length of $Q$. 
Observe that $r(Q)\sim\diam(B(Q))\sim\ell(Q)$.
Also, we call $z_Q$ the center of $Q$, and the cube $Q'\in \DD_{k-1}$ such that $Q'\supset Q$ the parent of $Q$.
 We set
$B_Q=28 \,B(Q)=B(z_Q,28\,r(Q))$.
%, so that 
%$W\cap \frac1{28}B_Q\subset Q\subset B_Q.$

We denote
$\DD^{db}=\bigcup_{k\geq0}\DD_k^{db}$.
Note that, in particular, from \rf{eqdob22} it follows that
\begin{equation}\label{eqdob*}
\omega^{p}(3B_{Q})\leq \omega^p(100B(Q))\leq C_0\,\omega^p(Q)\qquad\mbox{if $Q\in\DD^{db}.$}
\end{equation}
For this reason we will refer to the cubes from $\DD^{db}$ doubling.  As shown in \cite[Lemma 5.28]{David-Mattila}, every cube $R\in\DD$ can be covered $\omega^p$-a.e.
by a family of doubling cubes.
%\begin{lemma}\label{lemcobdob}
%Let $R\in\DD$. Suppose that the constants $A_0$ and $C_0$ in Lemma \ref{lemcubs} are
%chosen suitably. Then there exists a family of
%doubling cubes $\{Q_i\}_{i\in I}\subset \DD^{db}$, with
%$Q_i\subset R$ for all $i$, such that their union covers $\omega^p$-almost all $R$.
%\end{lemma}

\section{Good   and bad   cells. The  estimate of $\RR$ on good doubling cells}
\label{good}

We introduce the notions of  bad and good David--Mattila cells. First we need the $n$-dimensional Frostman measure.
Let $g\in L^1(\omega^p)$ be such that
$\omega^p|_E = g\,\HH^n|_{\partial\Omega}.$
Given $M>0$, let 
$$E_M= \{x\in\partial\Omega:M^{-1}\leq g(x)\leq M\}.$$
Take $M$ big enough so that $\omega^p(E_M)\geq \omega^p(E)/2$, say.
Consider an arbitrary compact set $F_M\subset E_M$ with $\omega^p(F_M)>0$. We will show that there exists $G_0\subset F_M$
with $\omega^p(G_0)>0$ which is $n$-rectifiable. Clearly, this suffices to prove that $\omega^p|_{E_M}$ is $n$-rectifiable,
and letting $M\to\infty$ we get the full $n$-rectifiability of $\omega^p|_E$.

Let $\mu$ be an $n$-dimensional Frostman measure for $F_M$. That is, $\mu$ is a non-zero Radon measure supported on $F_M$
such that 
$\mu(B(x,r))\leq C\,r^n$ for all $x\in\R^{n+1}$.
Further, by renormalizing $\mu$, we can assume that $\|\mu\|=1$. Of course the constant $C$ above will depend on 
$\HH^n_\infty(F_M)$, and the same may happen for all the constants $C$ to appear,  but this will not bother us. Notice that $\mu\ll\HH^n|_{F_M}\ll \omega^p$. In fact, for any set $H\subset F_M$,
\begin{equation}
\label{Frostman}
\mu(H)\leq C\,\HH^n_\infty(H)\leq C\,\HH^n(H)\leq C\,M\,\omega^p(H).
\end{equation}

The cell $Q\in \DD$ is called bad if it is a   maximal   cube satisfying one of the conditions below: 
\begin{itemize}
\item high density (HD): $\omega^p(3B_Q) \ge A \ell(Q)^n$, 
\,where  $A$\,\, is a suitably large number,
\item or low $\mu$ (LM): $\mu(B(Q)) \le \tau \omega^p(Q)$, where $\tau$ is a suitably small number.
\end{itemize}
  Any cube which is not contained in a bad cell will be called good. 

Notice that 
$ 
\sum_{  Q - bad,  \,Q\in\LM} \mu(Q) \leq \tau \sum_{  Q - bad,  \,Q\in\LM} \omega^p(Q) \leq \tau\,\|\omega\|=\tau=\tau\,\mu(F_M).$
Therefore, taking into account that $\tau$ is small, that $\mu$ is dominated by $\HH^n$ and that $\omega^p|_{F_M}$ and $\HH^n|_{F_M}$ are boundedly equivalent, we conclude that $\omega^p(F_M\setminus\bigcup_{  Q - bad,  \,Q\in LM} Q)>0$. Also,  taking into account that upper density of $\omega^p$, $\theta^*_{\omega^p}(x)=\limsup_{r\to 0} \frac{\omega^p(B(x,r)}{r^n},$ is finite
$\omega^p|_E$-a.e. (again because $\omega^p|_E\ll\HH^n|_E$), for $A$ large enough,   we can dispose of (HD) cubes as well. Ultimately, a little more careful consideration shows that $\omega^p(F_M\cap \bigcup_{Q\in \DD_0^{db}}Q\setminus\bigcup_{Q -bad} Q)>0$\footnote{Here we start working with a slightly enlarged collection of doubling cubes, but let us disregard these details in the present brief sketch.}.   
Notice also that the (HD) condition on bad cubes implies that 
\begin{equation}
\label{max}
  \M^n\omega^p(x)\lesssim A\quad \mbox{ for $\omega^p$-a.e.\ $x\in F_M\setminus \bigcup_{Q - bad} Q$.}  
\end{equation}
  It remains to concentrate on good cells contained in some cube from $\DD_0^{db}$. We start with 
\begin{lemma}
\label{main}
For any good doubling cell $R$ one has the estimate
$$
\left|\RR_{r(B_R)} \omega^p(x)\right| \lesssim C(A, M, \tau, {\rm dist}(p, \pom))\,\quad \mbox{ for any } x\in R.
$$
\end{lemma}

  The idea is to reduce the desired bound to certain estimates on the Green function of $\Omega$, which can be written as 
$$G(x,y) = \mathcal{E}(x-y) - \int_{\partial\Omega} \mathcal{E}(x-z)\,d\omega^y(z)\quad \mbox{ for $m$-a.e. $x\in\R^{n+1}$,}$$
where $\mathcal{E}$ is the fundamental solution for the Laplacian. Since the kernel of the Riesz transform is $K(x) = c_n\,\nabla \EE(x),$ essentially differentiating the expression above and using a trivial estimate on $\nabla \EE(x)$ in terms of distance from $p$ to $\pom$, we are left with considering the gradient of the Green function\footnote{In a careful argument, one has to utilize a smooth truncation of the kernel of the Riesz transform. The difference between the latter and $\RR_{r(B_R)} \omega^p(x^R)$ is fairly directly controlled by the $\M^n\omega^p$ which is in turn controlled by \eqref{max}.}. In fact, all we need is the bound $\frac1r\,| G(y,p)|\lesssim 1$ for all $y\in  B(x,r)\cap\Omega$, $r\approx r(B_Q)$.

  To this end, take $\delta>0$ from Bourgain's Lemma  
and fix a point $x^{R}$ such that $\mu(B(x^{R}, \delta r(R)))$ is maximal (or almost maximal) possible. Denote $B^{R}:=B(x^{R}, \delta r(R))$.  Then automatically,
$
\delta^n \mu(R)\lesssim \mu(B^{R})
$.
Thus, by Bourgain's lemma and by (\ref{Frostman}) we have
$
\inf_{z\in \partial (2B^{R})}  \omega^z(3B_{R}) \ge c\frac{\mu(R)}{r(R)^n}\,.
$
This very easily implies (for $n\ge 2$ at least) that for every $z\in \partial (2B^{R})$ and every $y\in B^{R}$ one has
$
\frac1{r(R)}G(z, y) \lesssim \frac{ \omega^z(3B_{R}) }{\mu(R)}
$. Now the maximal principle yields (we also use that $R$ is good and doubling in the last inequality)
\begin{equation}
\label{key}
\forall  y\in B^{R},\,\,\,\,\,\,\frac1{r(R)} G(y, p) \lesssim \frac{ \omega^p(3B_{R}) }{\mu(R)}\lesssim \tau^{-1}\,,
\end{equation}
as desired. The similar inequality can be proved for $n=1$ but requires more work.

  Consider now the case when $Q$ is good and non-doubling, $Q\notin \DD^{db}$, such that, in addition, $Q\subset R_0$, where $R_0\in \DD^{db}_0$. Let $R\supset Q$ be the cube from $\DD^{db}$ with minimal side length containing $Q$. Then for all $y\in B_Q$ we have
$$
|\RR_{r(B_Q)}\omega^p (y)| \leq |\RR_{r(B_{R})}\omega^p (y)| + C\,\sum_{P\in\DD: Q\subset P\subset R}\theta_{\omega^p}(2B_P),
$$
where $\theta_{\omega^p}(B):= \frac{\omega^p(B)}{r(B)^n}$. The first term is bounded by some constant as above since $R$ is good and doubling. The sum is bounded by $\theta_{\omega^p}(4B_{R})$ as for all such cubes 
$$
\omega^p(100B(P))\leq A_0^{-10n(J(P)-J(R)-1)}\omega^p(100B(R))
$$ 
see  \cite[Lemma 5.31]{David-Mattila}. Finally, $\theta_{\omega^p}(4B_{R})$ is bounded by $A$ as $R$ is not high density -- see (HD). This combines with Lemma~\ref{main} to yield the key theorem:
\begin{theorem}
\label{keyth} Let $G_M:= F_M\cap \bigcup_{Q\in \DD_0^{db}}Q\setminus\bigcup_{Q -bad} Q$. Then $\omega^p(G_M)>0$ and 
$$
\RR_*\omega^p(x) \lesssim C(A, M, \tau, {\rm dist}(p, \pom))\, \quad \mbox{ for  $\omega^p$-a.e. }x\in G_M. 
$$
\end{theorem}

\section{Harmonic measure: its singular integral and its rectifiability}
\label{endproof}

\begin{theorem} 
 \label{teo**}
Let $\sigma$ be a Radon measure with compact support on $\R^{n+1}$   and consider a $\sigma$-measura\-ble set
$G$ with $\sigma(G)>0$ such that
$$G\subset\{x\in \R^d: 
\M^n\sigma(x) < \infty \mbox{ and } \,\RR_* \sigma(x) <\infty\}.$$
Then there exists a Borel subset $G_0\subset G$ with $\sigma(G_0)>0$  such
that $\sup_{x\in G_0}\M^n\sigma|_{G_0}(x)<\infty$ and  $\RR_{\sigma|_{G_0}}$ is bounded in $L^2(\sigma|_{G_0})$.
\end{theorem}

This result follows from the  non-homogeneous $Tb$ theorem of Nazarov, Treil and Volberg in \cite{NTV} (see also \cite{Volberg}) in combination with the methods in \cite{Tolsa-pams}. For the detailed proof in the case of the Cauchy
transform, see \cite[Theorem 8.13]{Tolsa-llibre}. The same arguments with very minor modifications work for the Riesz transform.   This result applies to quite general antisymmetric Calderon-Zygmund operators.  
However, the next theorem uses very prominently that we work with the Riesz transform, and that its singularity is precisely $n$. It is the main result in \cite{NToV-pubmat}, which is based in its turn on the solution of (co-dimension $1$) David--Semmes conjecture in \cite{NToV}
%and on \cite{ENV} 
(see also \cite{HMM} in the context of uniform domains).

\begin{theorem}
\label{ntov}
Let $\sigma$ be Radon measure on $\R^{n+1}$ and $\sigma\ll\HH^n$. If  the vector $n$-dimensional Riesz transform $\RR $ is a bounded operator on $L^2(\sigma)$, then $\sigma$ is rectifiable.
\end{theorem}

By the latter theorem (or the David--L\'eger theorem  \cite{Leger} for $n=1$), we deduce
that $\omega^p|_{G_0}$ is $n$-rectifiable.

\end{document}